\documentclass[12pt]{extarticle}

\usepackage[margin=1.5in]{geometry}  
\usepackage{graphicx}              
\usepackage{amsmath}               
\usepackage{amsfonts}              
\usepackage{amsthm}                
\usepackage{framed}
\usepackage{url}

\newtheorem{thm}{Theorem}
\newtheorem{lem}[thm]{Lemma}

\begin{document}

\nocite{*}

\title{\bf Some Notes on Digit Strings in the Primes}

\author{\textsc{Adrian W. Dudek} \\ 
Mathematical Sciences Institute \\
The Australian National University \\ 
\texttt{adrian.dudek@anu.edu.au}}
\date{}

\maketitle

\begin{abstract}
\noindent Let $S$ be a string of $l$ decimal digits. We give an explicit upper bound on some prime $p$ whose decimal representation contains the string $S$. We also show, as a corollary of the Green-Tao theorem, that there are arbitrarily long arithmetic progressions of prime numbers all of whose decimal representations contain $S$. 
\end{abstract}

\section{Introduction}

It is an elementary exercise to show that if a series converges, then the sequence formed by its individual terms must tend to zero. The converse, however, does not follow; the archetype here is the harmonic series
\begin{equation} \label{harmonicseries}
1+\frac{1}{2} + \frac{1}{3} + \frac{1}{4} + \cdots
\end{equation}
which diverges to infinity. In 1914, Kempner \cite{kempner} proved that if one omits from the above series all terms whose denominator contains the digit 9, the resulting series will converge to a finite value. Kempner went on to show that this value was less than 80, and Baillie \cite{baillie} later gave the result to 20 decimal places; the sum is approximately $22.92067$. 

There is nothing special about the digit 9. One can, in fact, fix any finite string of digits and omit all terms from (\ref{harmonicseries}) whose denominator contains this string to get a convergent series. A proof of this can be found in Section 9.9 of Hardy and Wright's well-known text \cite{hardywright}. It should also be noted that Schmelzer and Baillie \cite{schmelzerbaillie} have provided an algorithm to compute the value of these series to high precision. It follows, from the convergence of such series, that we have the somewhat paradoxical statement that almost all integers contain any given string of digits. This is meant in the following way: if we fix a string $S$ and let $R(x)$ count the positive integers up to $x$ which do not contain the string $S$, then
$$\lim_{x \rightarrow \infty} \frac{R(x)}{x} = 0.$$
The purpose of this paper is to consider the consequences of this result in the setting of prime numbers. It is known that the series of the reciprocals of the primes diverges and, as this is a subseries of the harmonic series, it follows that there are infinitely many primes which contain any given string (see Behforooz \cite{behforooz}). Our first result gives an explicit upper bound on the least prime containing some fixed string. 

\begin{thm} \label{primebound}
Let $S$ be a string of $l$ digits. Then there exists a prime number $p$ containing $S$ within its decimal representation and satisfying the bound
\begin{equation} \label{thmeqn}
\frac{\log p}{\log \log p} \leq 5.7 l^2 \times 10^l.
\end{equation}
\end{thm}

Actually, the proof of this theorem gives something stronger, for it does not rely on the composition of the string but only its length. As such, for \textit{any} string $S$ of length $l$ there exists a prime $p \leq N$ such that $p$ contains $S$ in its decimal representation and 
\begin{equation} \label{forn}
\frac{\log N}{\log \log N} = 5.7 l^2 \times 10^l.
\end{equation}
This bound is considerably weak when contrasted with some calculations. Using text files containing the first two million primes (see Caldwell's tables \cite{caldwell}), we wrote a program in Java which, given a positive integer $l$, would return the least number $M$ such that all strings of length $l$ were contained in primes at most $M$. Table 1 compares these computations with the approximate values for $N$ we get by solving (\ref{forn}).

\begin{table}[ht] 
\caption{Experimentally determined values for $M$ against values for $\log N$ given by Theorem 1.} 
\label{tab:pixs} 
\centering 
\begin{tabular}{c c c } 
\\ \hline\hline 
$l$ & $M$ & $\log N$ \\[0.5 ex] \hline
$1$ & 83 & 330.7 \\
$2$ & 1847 & 22887.4 \\
$3$ & 50411 & 689676\\
$4$ & 793343 & $1.51 \cdot 10^7$ \\
$5$ & 9810001 & $2.77 \cdot 10^8$ \\
\hline\hline 
\end{tabular} 
\end{table}

Theorem \ref{primebound} thus seems an extremely weak result. There are some potential reasons for this; we see no grounds why certain strings should be favoured by primes nor others by composites, except for at the basic level where even-integer strings should not appear at the end of primes. There are also well-known divisibility rules which depend on decimal digits, though we see no way in which one would implement these. 

To this end, the proof of Theorem \ref{primebound} relied on a simple argument. If we consider the integers in the interval $[1,x]$ and let $\pi(x)$ denote the number of primes in this interval and $R(x)$ denote those integers which avoid some string $S$, then we can solve the inequality
$$\pi(x) > R(x)$$
using known estimates for $\pi(x)$. 

We can speculate somewhat on the correct order of Theorem \ref{primebound}, using an entertainingly loose argument based on the solution to the coupon collector's problem (see, for example, Isaac \cite{prob}): 

\begin{center}
\textit{Suppose that there is an urn of n different coupons, from which coupons are being collected, equally likely, with replacement. What is the expected number of collections required to guarantee that all different coupons have been collected?}
\end{center}

It is known that the expected number is asymptotic to $n \log n$. Our speculation on the correct order of Theorem \ref{primebound} will be questionable for we suppose that, for a given prime with at least $l$ digits, it is equally likely that it will contain any string of length $l$ i.e. there is no bias amongst strings. We must also suppose that there is at most one string in each prime, which is risible to say the least.

Then, consider that for each prime following $10^{l-1}$, we wish to draw from it a string of length $l$. There are $9 \cdot 10^{l-1}$ such strings, and so we expect to have to go through the first
$$9 \cdot 10^{l-1} \log(9 \cdot 10^{l-1})$$
primes which exceed $10^{l-1}$ in order to collect all strings of the given length. By the prime number theorem, there are approximately
$$\frac{10^{l-1}}{(l-1) \log 10}$$
primes which do not exceed $10^{l-1}$. Thus, we should have drawn all strings of length $l$ by $N$, where
$$\pi(N) \approx \frac{10^{l-1}}{(l-1) \log 10} + 9 \cdot 10^{l-1} \log(9 \cdot 10^{l-1}).$$
It is straightforward enough to take $l$ to be large so that we may neglect some terms and then invert the formula to get that
$$N \sim k l^2 10^l$$
for some constant $k$. One may compare this with (\ref{thmeqn}). Though there is not an ample amount of computational data, the values for $M$ in Table 1 do not strongly contradict this.

In Section 3, we give a short proof for the following theorem, which turns out to be a corollary of Theorem 1.2 of Green and Tao's \cite{greentao} celebrated paper and our proof of Theorem \ref{primebound}.

\begin{thm} \label{greentaocorollary}
Let $S$ be a string of length $l$. Then there are arbitrarily long arithmetic progressions of prime numbers whose decimals contain the string $S$.
\end{thm}

\section{Proof of Theorem \ref{primebound}}

Following Hardy and Wright \cite{hardywright}, we will consider integers with decimal representations in base $r$. This will allow for simplicity, for a string of length $l$ will correspond to a single digit if we choose $r = 10^l$. We denote by $v$ an integer whose decimal representation avoids the digit $b$. Then the number of such integers for which $r^{l-1} \leq v < r^l$ is $(r-1)^{l}$ in the case that $b=0$ and $(r-2)(r-1)^{l-1}$ if $b \neq 0$. In either case, the count does not exceed $(r-1)^l$. Thus, if 
$$r^{k-1} \leq x < r^k,$$
then the number $R(x)$ of $v$ up to $x$ does not exceed
$$(r-1)+(r-1)^2 + \cdots + (r-1)^k = \frac{(r-1)^{k+1}-(r-1)}{r-2} < \frac{(r-1)^{k+1}}{r-2}.$$
Now, when there are more primes than there are numbers whose decimal representations avoid $b$, we have that there must be some prime whose representation contains $b$. Therefore, we wish to solve the inequality
$$\pi(x) > R(x)$$
for $x$, where $\pi(x)$ denotes the number of primes not exceeding $x$. We can do this by continuously strengthening the inequality to get a simpler inequality involving $x$. We can employ Corollary 1 of Rosser and Schoenfeld \cite{rosserschoenfeld1962}, namely that
\begin{equation} \label{piestimate}
\pi(x) > \frac{x}{\log x}
\end{equation}
for all $ x \geq 17$, with our upper bound for $R(x)$, so that we may instead solve the stronger inequality
\begin{equation} \label{one}
\frac{x}{\log x} > \frac{(r-1)^{k+1}}{r-2}
\end{equation}
for $x \geq 17$. The reader should note that there are sharper estimates than (\ref{piestimate}) available (see Dusart \cite{dusart} for instance) but these do not lead to a significant improvement of the main result; an improvement would only be made on the error terms which we have chosen to neglect. 

Using the fact that $r^{k-1} \leq x$, we can rearrange (\ref{one}) to get

$$\log x < \frac{r-2}{r(r-1)} \Big( \frac{r}{r-1} \Big)^k.$$
Taking the logarithm of both sides and using 

$$\frac{\log x}{\log r} < k$$
we can get

$$\log \log x < \log\Big( \frac{r-2}{r(r-1)} \Big) + \log x \Big( 1- \frac{\log(r-1)}{\log r} \Big).$$
It remains to rearrange and use the bound

$$\log r - \log(r-1) < \frac{1}{r}$$
along with the trivial bound $\log \log x > 1$ when $x > \exp(e)$ to get that

$$\frac{\log x}{\log \log x} > r \log^2 r \Big(1+\frac{1 + \log(\frac{r-1}{r-2})}{\log r}\Big).$$
The above inequality implies that $\pi(x) > R(x)$, and so there must exist some prime $p$ containing the digit $b$ which satisfies

$$\frac{\log p}{\log \log p} \leq r \log^2 r \Big(1+\frac{1 + \log(\frac{r-1}{r-2})}{\log r}\Big).$$

Now note that a string of length $l$ corresponds to a single digit in base $r=10^l$, so we have that for any string of length $l$ there exists a prime $p$ whose decimal contains the string and satisfies 

$$\frac{\log p}{\log \log p} \leq 5.7 l^2 \times 10^l$$
so long as $l \geq 6$. Table 1 confirms this bound for the remaining values of $l$, and this completes the proof of Theorem \ref{primebound}.

\section{Proof of Theorem \ref{greentaocorollary}}

We require the following lemma, which is Theorem 1.2 of Green and Tao \cite{greentao}.

\begin{lem}
Let $A$ be any subset of the prime numbers of positive relative upper density, thus $\limsup_{N \rightarrow \infty} \pi(N)^{-1} | A \cap [1,N]|>0$. Then $A$ contains infinitely many arithmetic progressions of length $k$ for all $k$.
\end{lem}

Theorem 2 clearly follows from showing that for any string $S$, the set $A$ of primes whose decimal contains the string $S$ has positive relative upper density. The set $A$ actually has a density of 1. We let $B$ denote the set of primes whose decimal does not contain the string $S$, and we write $A(N) = A \cap [1,N]$ and $B(N) =  B \cap [1,N].$

From the proof of Theorem \ref{primebound}, we have that

$$\frac{R(N)}{N} < \frac{r (r-1)}{r-2} \Big( \frac{r-1}{r} \Big)^k.$$
for $r^{k-1} \leq N < r^k$. Clearly,

\begin{eqnarray*}
\frac{ | A \cap [1,N]|}{\pi(N)} & = & 1 - \frac{|B(N)|}{\pi(N)} \\
& = & 1 + O\bigg( \frac{R(N)}{\pi(N)}  \bigg)\\
& = & 1 + O\bigg( \frac{R(N) \log N}{N} \bigg).
\end{eqnarray*}
It is straightforward to check that
$$R(N) = o(N / \log N)$$
as $N \rightarrow \infty$ and so the error term tends to zero. This completes the proof of Theorem \ref{greentaocorollary}.

\newpage

\bibliographystyle{plain}

\bibliography{biblio}

\end{document}